# On the Combinatorial Inverse Monoid $IO_3$

**Dr. Emil Daniel Schwab, Professor,**
**Romero Efren, MS Student,**
University of Texas at El Paso, U.S.A.

**ABSTRACT.** In this paper we compute the Mobius category and the Mobius function μ of the combinatorial inverse monoid $IO_3$ of all order preserving partial bijections on the set $M_3=\{1,2,3\}$. This category is the reduced standard division category $C_F(IO_3)$ relative to an idempotent transversal F of the $D$-classes of $IO_3$. It is a special Mobius category namely a Mobius category in which the values of $\mu(\alpha)$ are 1 and –1 for any morphism α of $C_F(IO_3)$.

## 1. Introduction

Let S be a monoid. Recall the definition of the Green's relations *R,L,H* and *D*. For all s,t∈S,

$$sRt \quad \text{if and only if} \quad sS = tS$$
$$sLt \quad \text{if and only if} \quad Ss = St,$$

and

$$H = R \cap L \quad ; \quad D = R \circ L = L \circ R$$

A semigroup (monoid) is said to be an inverse semigroup (monoid) if, for every element s, there is an unique t (called the inverse of s) such that s=sts and t=tst. The unique inverse of s shall be denoted $s^{-1}$. The following are two equivalent definitions of an inverse monoid S:
(i) every element s in S has at least one inverse and the idempotents of S commute;
(ii) each *R*-class and each *L*-class of S has exactle one idempotent.

213

Notice that $ss^{-1}$ and $s^{-1}s$ are idempotents for all $s \in S$. As usual, $E(S)$ denote the set of all idempotents of S.

For any inverse monoid S a relation may be defined which makes S a poset. The relation is given by:

$$s \leq t \quad \text{if and only if} \quad s \in E(S)t \quad \text{for any } s,t \in S$$

The following conditions are equivalent to the above partial order (called the natural partial order on S):
1. $s = et$ where $e \in E(S)$
2. $ss^{-1} = st^{-1}$
3. $s^{-1}s = t^{-1}s$
4. $s = st^{-1}s$.

Applying this relation to the commutative semigroup of idempotents $E(S)$, the relation becomes:

$$e \leq f \quad \text{if and only if} \quad e = ef = fe \quad \text{for any } e, f \in E(S),$$

which is a semilattice. The Green's equivalence relations $R, L$ and $H$ on the inverse monoid S are given by:

$$sRt \quad \text{if and only if} \quad ss^{-1} = tt^{-1}$$
$$sLt \quad \text{if and only if} \quad s^{-1}s = t^{-1}t$$
$$sHt \quad \text{if and only if} \quad ss^{-1} = tt^{-1} \text{ and } s^{-1}s = t^{-1}t$$

and two idempotents e and f of the inverse monoid S are *D*-related,

$eDf$, if and only if there is an element $s \in S$ such that $e = s^{-1}s$ and $f = ss^{-1}$.

By combinatorial inverse monoids we mean inverse monoids whose al subgroups are trivial. An inverse monoid is combinatorial if and only if *H* is the equality relation.

The set of all partial bijections from a set M to itself forms an inverse monoid *I*(M) called the symmetric inverse monoid of M (the elements of *I*(M) are bijections between subsets of M and the product of two such bijections is their product as partial transformations). An element f of *I*(M) is idempotent if and only if $f = 1_A$ for some subset A of M. The semilattice of idempotents of *I*(M) is isomorphic to the semilattice of subsets of M under intersection. The set of all elements of *I*(M) having domain and range equal to M is a subgroup of $I$(M) called the symmetric group on the set M. Thus $I$(M) ($|M|>1$) is not a combinatorial inverse monoid. But the subsemigroup $IO_n$ of order preserving partial bijections of a chain with n elements is combinatorial.

Leech [Lee87] showed that special categories (called division categories) could be used to generalize Clifford's constructions of bisimple inverse monoid from a class of right cancellative monoids by replacing the monoid by a suitable category. In fact inverse monoids and division categories are equivalent structures: each can be reconstructed from the other. We show how to construct a division category from an inverse monoid . But first we give the definition of a division category.

A division category is a small category C with the following properties:
(1) every morphism of C is an epimorphism;
(2) C has an object I (called a quasi initial object) such that Hom(I,X)≠Φ for any X∈ObC;
(3) C has pushouts.

If S is an inverse monoid, the category C(S) defined by:
- ObC(S)=E(S);
- Hom(e,f)={(s,e)| s∈S such that $s^{-1}s \leq e$ and $ss^{-1} = f$ }
- (t,f)·(s,e)=(ts,e) is the composition of two morphisms (s,e):e→f and (t,f):f→g

is a division category (with 1 a quasi initial object) called the standard division category of S. Now, if F is an idempotent transversal of the *D*-classes of S such that 1∈F, then the full subcategory $C_F(S)$ of C(S) given by $ObC_F(S)=F$ is another division category called a reduced standard division category of S. For two idempotent transversals F and F' (with 1∈F and 1∈F') , $C_F(S)$ and $C_{F'}(S)$ are isomorphic (see [JL99]). The first author [Sch04a] showed that a reduced standard division category $C_F(S)$ of an inverse monoid S is a Mobius category if and only if S is combinatorial and $(E(S),\leq)$ is locally finite. Mobius categories were introduced by Leroux [Ler75]. They were created as a general program to extend the theory of Mobius functions. We refer the reader to [CLL80],[Sch04a] and [Sch04b] for all standard definitions and results from Mobius category theory. Given a combinatorial inverse monoid S with Mobius reduced standard division category $C_F(S)$, the Mobius function $\mu:MorC_F(S) \to \mathbf{C}$ is defined by

$$\mu(s,e) = \mu_{E(eSe)}([s^{-1}s,e]_{E(eSe)}),$$

where $\mu_{E(eSe)}$ is the Mobius function of the locally finite lattice $(E(eSe),\leq)$ and $[s^{-1}s,e]_{E(eSe)}$ is the interval $\{f \in E(eSe) \mid s^{-1}s \leq f \leq e\}$.

215

## 2. The reduced standard division category of $IO_3$

In order to illustrate the construction of a standard division category $C_F(S)$ from an inverse monoid S consider the set $M_3 = \{1,2,3\}$. Then for the symmetric inverse monoid $I(M_3)$ we have $|I(M_3)| = 34$. The notation for the partial bijection $f \in I(M_3)$ where f(2)=1, f(3)=3 and 1 has no image will be

$$\begin{pmatrix} 1 & 2 & 3 \\ - & 1 & 3 \end{pmatrix}$$

(where the − indicates that 1 has no image). The notation will make for easy computation of compositions (identical to that of composition of permutation groups), however in order to facilitate the display of tables of composition an abbreviated notation will be adopted.

This notation for the above example will be [21]. This indicates that 2 maps into 1, 1 has no image, and 3 maps into 3. To make the criteria for the proper notation of and element $f \in I(M_3)$ more concrete note the following conventions:
1. [abc] indicates that f(a)=b and f(b)=c;
2. The last digit within a set of brackets has no image;
3. If the image of a is a then omit it.

Following these guidelines the partial bijections

$$\begin{pmatrix} a & b & c \\ b & c & - \end{pmatrix}, \begin{pmatrix} a & b & c \\ a & c & - \end{pmatrix}, \begin{pmatrix} a & b & c \\ a & - & b \end{pmatrix}$$

will be written as [abc], [b,c] and [cb] respectively. The only partial bijections that will not adopt this convention will be the zero and identity partial bijections which will be written as 0 and i.

Now, consider the set

$IO_3 = \{f \in I(M_3) \mid f \text{ is order preserving partial bijection}\}$.

The elements of $IO_3$ are

$IO_3 = \{[1],[2],[3],[12],[21],[23],[32],[123],[321],[1][2],[1][3],[2][3],$
$\quad [31][2],[13][2],[32][1],[23][1],[12][3],[21][3],0,i\}$.

Tables 2.1-2.4 list the composition of partial bijections in $IO_3$. Table 2.5 lists the inverse of each element of $IO_3$. The only elements not represented in tables 2.1-2.4 are 0 and i as their composition with any other element of $IO_3$ is clear.

*Table 2.1*

|        | [2][3] | [12][3] | [13][2] | [21][3] | [1][3] | [23][1] | [31][2] | [32][1] | [1][2] |
|--------|--------|---------|---------|---------|--------|---------|---------|---------|--------|
| [2][3] | [2][3] | 0       | 0       | [21][3] | 0      | 0       | [31][2] | 0       | 0      |
| [12][3]| [12][3]| 0       | 0       | [1][3]  | 0      | 0       | [32][1] | 0       | 0      |
| [13][2]| [13][2]| 0       | 0       | [23][1] | 0      | 0       | [1][2]  | 0       | 0      |
| [21][3]| 0      | [2][3]  | 0       | 0       | [21][3]| 0       | 0       | [31][2] | 0      |
| [1][3] | 0      | [12][3] | 0       | 0       | [1][3] | 0       | 0       | [32][1] | 0      |
| [23][1]| 0      | [13][2] | 0       | 0       | [23][1]| 0       | 0       | [1][2]  | 0      |
| [31][2]| 0      | 0       | [2][3]  | 0       | 0      | [21][3] | 0       | 0       | [31][2]|
| [32][1]| 0      | 0       | [12][3] | 0       | 0      | [1][3]  | 0       | 0       | [32][1]|
| [1][2] | 0      | 0       | [13][2] | 0       | 0      | [23][1] | 0       | 0       | [1][2] |

*Table 2.2*

|        | [3]    | [23]   | [123]  | [32]   | [2]    | [12]   | [321]  | [21]   | [1]    |
|--------|--------|--------|--------|--------|--------|--------|--------|--------|--------|
| [2][3] | [2][3] | [2][3] | 0      | [2][3] | [2][3] | 0      | [21][3]| [21][3]| 0      |
| [12][3]| [12][3]| [12][3]| 0      | [12][3]| [12][3]| 0      | [1][3] | [1][3] | 0      |
| [13][2]| [13][2]| [13][2]| 0      | [13][2]| [13][2]| 0      | [23][1]| [23][1]| 0      |
| [21][3]| [21][3]| 0      | [2][3] | [31][2]| 0      | [2][3] | [31][2]| 0      | [21][3]|
| [1][3] | [1][3] | 0      | [12][3]| [32][1]| 0      | [12][3]| [32][1]| 0      | [1][3] |
| [23][1]| [23][1]| 0      | [13][2]| [1][2] | 0      | [13][2]| [1][2] | 0      | [23][1]|
| [31][2]| 0      | [21][3]| [21][3]| 0      | [31][2]| [31][2]| 0      | [31][2]| [31][2]|
| [32][1]| 0      | [1][3] | [1][3] | 0      | [32][1]| [32][1]| 0      | [32][1]| [32][1]|
| [1][2] | 0      | [23][1]| [23][1]| 0      | [1][2] | [1][2] | 0      | [1][2] | [1][2] |

*Table 2.3*

|        | [2][3] | [12][3] | [13][2] | [21][3] | [1][3] | [23][1] | [31][2] | [32][1] | [1][2] |
|--------|--------|---------|---------|---------|--------|---------|---------|---------|--------|
| [3]    | [2][3] | [12][3] | 0       | [21][3] | [1][3] | 0       | [31][2] | [32][1] | 0      |
| [23]   | [2][3] | [13][2] | 0       | [21][3] | [23][1]| 0       | [31][2] | [1][2]  | 0      |
| [123]  | [12][3]| [13][2] | 0       | [1][3]  | [23][1]| 0       | [32][1] | [1][2]  | 0      |
| [32]   | [2][3] | 0       | [12][3] | [21][3] | 0      | [1][3]  | [31][2] | 0       | [32][1]|
| [2]    | [2][3] | 0       | [13][2] | [21][3] | 0      | [23][1] | [31][2] | 0       | [1][2] |
| [12]   | [12][3]| 0       | [13][2] | [1][3]  | 0      | [23][1] | [32][1] | 0       | [1][2] |
| [321]  | 0      | [2][3]  | [12][3] | 0       | [21][3]| [1][3]  | 0       | [31][2] | [32][1]|
| [21]   | 0      | [2][3]  | [13][2] | 0       | [21][3]| [23][1] | 0       | [31][2] | [1][2] |
| [1]    | 0      | [12][3] | [13][2] | 0       | [1][3] | [23][1] | 0       | [32][1] | [1][2] |

*Table 2.4*

|       | [3]    | [23]   | [123]  | [32]   | [2]    | [12]   | [321]  | [21]   | [1]    |
|-------|--------|--------|--------|--------|--------|--------|--------|--------|--------|
| [3]   | [3]    | [2][3] | [12][3]| [32]   | [2][3] | [12][3]| [321]  | [21][3]| [1][3] |
| [23]  | [23]   | [2][3] | [13][2]| [2]    | [2][3] | [13][2]| [21]   | [21][3]| [23][1]|
| [123] | [123]  | [12][3]| [13][2]| [12]   | [12][3]| [13][2]| [1]    | [1][3] | [23][1]|
| [32]  | [2][3] | [3]    | [1][3] | [2][3] | [32]   | [32][1]| [21][3]| [321]  | [32][1]|
| [2]   | [2][3] | [23]   | [23][1]| [2][3] | [2]    | [1][2] | [21][3]| [21]   | [1][2] |
| [12]  | [12][3]| [123]  | [23][1]| [12][3]| [12]   | [1][2] | [1][3] | [1]    | [1][2] |
| [321] | [21][3]| [1][3] | [3]    | [31][2]| [32][1]| [32]   | [31][2]| [32][1]| [321]  |
| [21]  | [21][3]| [23][1]| [23]   | [31][2]| [1][2] | [2]    | [31][2]| [1][2] | [21]   |
| [1]   | [1][3] | [23][1]| [123]  | [32][1]| [1][2] | [12]   | [32][1]| [1][2] | [1]    |

*Table 2.5*

| s     | $s^{-1}$ | s      | $s^{-1}$ |
|-------|----------|--------|----------|
| [1]   | [1]      | [2][3] | [2][3]   |
| [2]   | [2]      | [1][2] | [1][2]   |
| [3]   | [3]      | [1][3] | [1][3]   |
| [12]  | [21]     | [13][2]| [31][2]  |
| [21]  | [12]     | [31][2]| [13][2]  |
| [23]  | [32]     | [12][3]| [21][3]  |
| [32]  | [23]     | [21][3]| [12][3]  |
| [123] | [321]    | [23][1]| [32][1]  |
| [321] | [123]    | [32][1]| [23][1]  |

The conclusion drawn from the tables above is that $IO_3$ is indeed an inverse monoid, and the idempotents of $IO_3$ are

$$E(S) = \{0, i, [1], [2], [3], [1][2], [1][3], [2][3]\}.$$

Now, the Cayley table of multiplication of idempotents is the following one:

***Table 2.6:*** *Cayley table of $E(IO_3)$*

|        | 0 | i      | [1]    | [2]    | [3]    | [1][2] | [1][3] | [2][3] |
|--------|---|--------|--------|--------|--------|--------|--------|--------|
| 0      | 0 | 0      | 0      | 0      | 0      | 0      | 0      | 0      |
| i      | 0 | i      | [1]    | [2]    | [3]    | [1][2] | [1][3] | [2][3] |
| [1]    | 0 | [1]    | [1]    | [1][2] | [1][3] | [1][2] | [1][3] | 0      |
| [2]    | 0 | [2]    | [1][2] | [2]    | [2][3] | [1][2] | 0      | [2][3] |
| [3]    | 0 | [3]    | [1][3] | [2][3] | [3]    | 0      | [1][3] | [2][3] |
| [1][2] | 0 | [1][2] | [1][2] | [1][2] | 0      | [1][2] | 0      | 0      |
| [1][3] | 0 | [1][3] | [1][3] | 0      | [1][3] | 0      | [1][3] | 0      |
| [2][3] | 0 | [2][3] | 0      | [2][3] | [2][3] | 0      | 0      | [2][3] |

Using the Cayley table 2.6, we can see that:

$$0 \leq 0, i, [1], [2], [3], [1][2], [1][3], [2][3]$$
$$i \leq i$$
$$[1] \leq i, [1]$$
$$[2] \leq i, [2]$$
$$[3] \leq i, [3]$$
$$[1][2] \leq i, [1], [2], [1][2]$$
$$[1][3] \leq i, [1], [3], [1][3]$$
$$[2][3] \leq i, [2], [3], [2][3].$$

The Hasse diagram of $(E(IO_3), \leq)$ is:

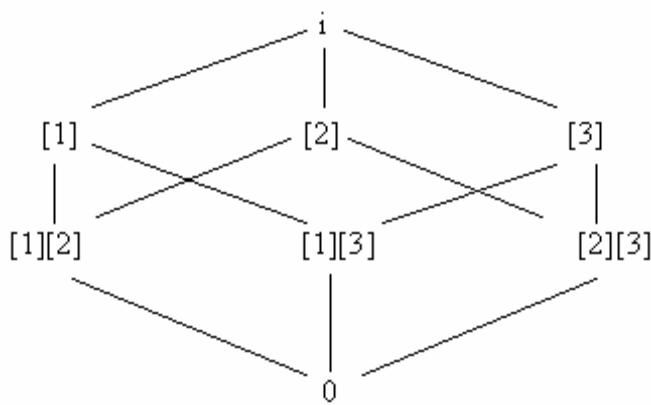

In order to construct the standard division category $C_F(IO_3)$ from $IO_3$ the products $s^{-1}s$ and $ss^{-1}$ must be known for all $s \in IO_3$. Table 2.7 provides these products.

*Table 2.7*

| s | $s^{-1}s$ | $ss^{-1}$ | s | $s^{-1}s$ | $ss^{-1}$ |
|---|---|---|---|---|---|
| [12] | [2] | [1] | [1][2] | [1][2] | [2][3] |
| [21] | [1] | [2] | [32][1] | [1][2] | [1][3] |
| [23] | [3] | [2] | [23][1] | [1][3] | [1][2] |
| [32] | [2] | [3] | [1][2] | [1][2] | [1][2] |
| [123] | [3] | [1] | [2][3] | [2][3] | [2][3] |
| [321] | [1] | [3] | [1][3] | [1][3] | [1][3] |
| [12][3] | [2][3] | [1][3] | [1] | [1] | [1] |
| [21][3] | [1][3] | [2][3] | [2] | [2] | [2] |
| [13][2] | [2][3] | [1][2] | [3] | [3] | [3] |

Since in an inverse semigroup two idempotents e and f are *D*-related if and only if there is an element s such that $e = s^{-1}s$ and $f = ss^{-1}$, we can see that [1],[2], and [3] are *D*-related and [1][2], [1][3] and [2][3] are *D*-related. Thus,

$$F = \{0, [1], [1][2], i\}$$

is an idempotent transversal of the *D*-classes of $IO_3$ with the identity i∈F.

The objects of the reduced standard division category $C_F(IO_3)$ are the elements of F: 0, [1], [1][2] and i. If e,f∈F, the set of morphisms from e to f is given by:

$$Hom(e, f) = \{(s,e) \in F \times S \mid s^{-1}s \leq e \text{ and } ss^{-1} = f\}.$$

It is straightforward to calculate the following table:

***Table 2.8:*** *Morphisms of $C_F(IO_3)$*

| s | $ss^{-1}$ | $s^{-1}s \leq e$ and $e \in F$ | $(s,e) \in Hom(e,f)$ |
|---|---|---|---|
| [12] | [1] | [2] $\leq$ i | ([12],i)$\in$Hom(i,[1]) |
| [21] | [2]$\notin$F | --- | --- |
| [23] | [2]$\notin$F | --- | --- |
| [32] | [3]$\notin$F | --- | --- |
| [123] | [1] | [3] $\leq$ i | ([123],i)$\in$Hom(i,[1]) |
| [321] | [3]$\notin$F | --- | --- |
| [12][3] | [1][3]$\notin$F | --- | --- |
| [21][3] | [2][3]$\notin$F | --- | --- |
| [13][2] | [1][2] | [2][3] $\leq$ i | ([13][2],i)$\in$Hom(i,[1][2]) |
| [31][2] | [2][3]$\notin$F | --- | --- |
| [32][1] | [1][3]$\notin$F | --- | --- |
| [23][1] | [1][2] | [1][3] $\leq$ [1][2], [1], i | ([23][1],[1])$\in$Hom([1],[1][2]) <br> ([23],i)$\in$Hom(i,[1][2]) |
| [1][2] | [1][2] | [1][2] $\leq$ [1][2], [1], i | ([1][2],[1][2])$\in$Hom([1][2],[1][2]) <br> ([1][2],[1])$\in$Hom([1],[1][2]) <br> ([1][2],i)$\in$Hom(i,[1][2]) |
| [2][3] | [2][3]$\notin$F | --- | --- |
| [1][3] | [1][3]$\notin$F | --- | --- |
| [1] | [1] | [1] $\leq$ [1], i | ([1],[1])$\in$Hom([1],[1]) <br> ([1],i)$\in$Hom(i,[1]) |
| [2] | [2]$\notin$F | --- | --- |
| [3] | [3]$\notin$F | --- | --- |
| 0 | 0 | 0 $\leq$ 0, [1][2], [1], i | (0,0)$\in$Hom(0,0) <br> (0,[1][2])$\in$Hom([1][2],0) <br> (0,[1])$\in$Hom([1],0) <br> (0,i)$\in$Hom(i,0) |
| i | i | i $\leq$ i | (i,i)$\in$Hom(i,i) |

It follows that:

$$Hom(i,[1][2]) = \{([13][2],i), ([23][1],i), ([1][2],i)\}$$
$$Hom(i,[1]) = \{([1],i), ([12],i), ([123],i)\}$$
$$Hom(i,0) = \{(0,i)\}$$
$$Hom([1],[1][2]) = \{([23][1],[1]), ([1][2],[1])\}$$
$$Hom([1],0) = \{(0,[1])\}$$
$$Hom([1][2],0) = \{(0,[1][2])\}$$

and

Hom(i,i) = {(i,i)}; Hom([1],[1]) = {([1],[1])}; Hom([1][2],[1][2]) = {([1][2],[1][2])}; Hom(0,0) = {(0,0)}.

Now, the reduced standard division category $C_F(IO_3)$ is described by the graph:

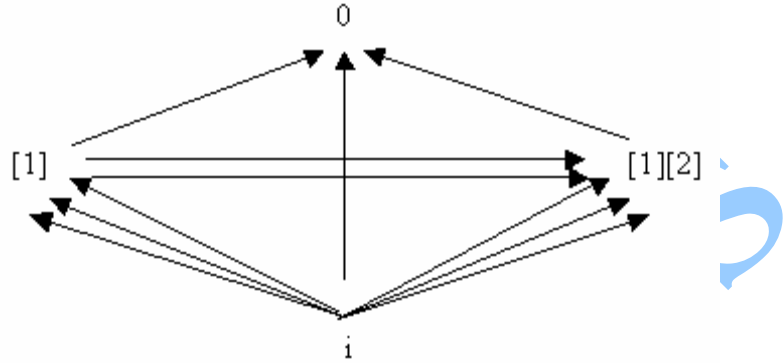

The morphism compositions are:

$$([23][1],[1]) \cdot ([1],i) = ([23][1],i)$$
$$([23][1],[1]) \cdot ([12],i) = ([13][2],i)$$
$$([23][1],[1]) \cdot ([123],i) = ([13][2],i)$$
$$([1][2],[1]) \cdot ([1],i) = ([1][2],i)$$
$$([1][2],[1]) \cdot ([12],i) = ([1][2],i)$$
$$([1][2],[1]) \cdot ([123],i) = ([23][1],i).$$

Any morphism of $C_F(IO_3)$ is an epimorphism and i is a quasi-initial object. But, not any morphism of $C_F(IO_3)$ is a monomorphism. For example the morphisms:

$$i \to 0, \quad [1] \to 0, \quad [1][2] \to 0$$

are not monomrphisms. Now, it is straightforward to check that for $IO_3$ the Green relation $H$ is the equality relation. So, $IO_3$ is combinatorial and therefore $C_F(IO_3)$ is a Mobius category.

## 3. The Mobius function

A decomposition finite category is a small category having finitely many decompositions α=βγ for any morphism α. The incidence algebra A(C) of a decomposition finite category C is the **C**-algebra of all complex valued functions $\xi: \text{Mor}C \to \mathbf{C}$ with the usual vector space structure and multiplication given by convolution:

$$(\xi * \eta)(\alpha) = \sum_{\beta\gamma=\alpha} \xi(\beta) \cdot \eta(\gamma).$$

The identity element of A(C) is $\delta$ given by:

$$\delta(\alpha) = \begin{cases} 1 & \text{if } \alpha \text{ is an identity morphism} \\ 0 & \text{otherwise} \end{cases}$$

A Mobius category is a decomposition finite category C such that an incidence function $\xi$ in A(C) has a convolution inverse if and only if $\xi(\alpha) \neq 0$ for any identity morphism $\alpha$. The Mobius function $\mu$ of a Mobius category C is the convolution inverse of the zeta function $\zeta$ (where $\zeta(\alpha)=1$ for any morphism $\alpha$ of C). The Mobius inversion formula is then nothing but the statement:

$$\eta = \xi * \zeta \iff \xi = \eta * \mu.$$

One way to find the Mobius function of $C_F(IO_3)$ is to apply the definition of Mobius function. For an identity morphism $\alpha$ of $C_F(IO_3)$ we have:

$$1 = \delta(\alpha) = (\zeta * \mu)(\alpha) = \zeta(\alpha) \cdot \mu(\alpha) = \mu(\alpha).$$

The following morphisms:

$$[1][2] \longrightarrow 0, \quad [1] \rightrightarrows [1][2], \quad i \rightrightarrows [1]$$

are non-identity indecomposable morphisms of $C_F(IO_3)$ (a morpfism $\beta$ is called indecomposable if $\beta=\beta'\beta''$ implies $\beta'$ is an identity or $\beta''$ is an identity). For a non-identity indecomposable morphism $\beta$ of $C_F(IO_3)$, we have:

$$0 = \delta(\beta) = (\zeta * \mu)(\beta) = \zeta(\beta) \cdot \mu(1) + \zeta(1) \cdot \mu(\beta) = 1 + \mu(\beta),$$

that is

$$\mu(\beta) = -1.$$

Now for the morphism $[1] \to 0$ we obtain:

$$0 = \delta([1] \to 0) = (\zeta * \mu)([1] \to 0) = \zeta([1] \to 0) \cdot \mu(1_{[1]}) + \zeta([1][2] \to 0) \cdot \mu(\beta') +$$
$$+ \zeta([1][2] \to 0) \cdot \mu(\beta'') + \zeta(1_0) \cdot \mu([1] \to 0),$$

where

$$\beta' = ([23][1],[1]) \quad and \quad \beta'' = ([1][2],[1]).$$

It follows:

$$0 = 1 - 1 - 1 + \mu([1] \to 0),$$

that is,

$$\mu([1] \to 0) = 1.$$

For any morphism $i \to [1][2]$ from i to $[1][2]$,

$$0 = \delta(i \to [1][2]) = (\zeta * \mu)(i \to [1][2]) = \zeta(i \to [1][2]) \cdot \mu(1_i) + \zeta([1] \to [1][2]) \cdot \mu(\beta') +$$
$$+ \zeta([1] \to [1][2]) \cdot \mu(\beta'') + \zeta(1_{[1][2]}) \cdot \mu(i \to [1][2])$$

where β' and β'' are two (non-identity indecomposable) morphisms from i to [1]. Thus,

$$0 = 1 - 1 - 1 + \mu(i \to [1][2]).$$

Therefore,

$$\mu(i \to [1][2]) = 1.$$

Finally,

$$0 = \delta(i \to 0) = (\zeta * \mu)(i \to 0) = \zeta(i \to 0) \cdot \mu(1_i) + \zeta([1] \to 0) \cdot \mu(i \xrightarrow{\beta} [1]) +$$
$$+ \zeta([1] \to 0) \cdot \mu(i \xrightarrow{\beta'} [1]) + \zeta([1] \to 0) \cdot \mu(i \xrightarrow{\beta''} [1]) + \zeta([1][2] \to 0) \cdot \mu(i \xrightarrow{\gamma} [1][2]) +$$
$$+ \zeta([1][2] \to 0) \cdot \mu(i \xrightarrow{\gamma'} [1][2]) + \zeta([1][2] \to 0) \cdot \mu(i \xrightarrow{\gamma''} [1][2]) + \zeta(1_0) \cdot \mu(i \to 0) =$$
$$= 1 - 1 - 1 - 1 + 1 + 1 + 1 + \mu(i \to 0).$$

Thus,

$$\mu(i \to 0) = -1.$$

Consequently the Mobius function of the Mobius-division category $C_F(IO_3)$ is given by:

$$\mu(\alpha) = \begin{cases} -1 & \text{if } \alpha = i \to 0 \text{ or } \alpha \text{ is a non-identity indecomposable morphism} \\ 1 & \text{otherwise} \end{cases}$$

**Remark.** *A Mobius category C is a special Mobius category if the values of $\mu(\alpha)$, as $\alpha$ ranges over the morhisms of C, are 1 and –1. Thus $C_F(IO_3)$ is a special Mobius category. The first author proved [Sch06] that the reduced standard division category $C_F(IO_n)$ is a special Mobius category for any positive integer n. This special Mobius category is isomorphic to the category $C_n$ of strictly increasing sequences of integers with indices:*
 - $ObC_n = \{0, 1, 2, \ldots, n\}$

 - $Hom(i, j) = \begin{cases} \{(a_1, a_2, \ldots, a_j)^i \mid a_1, a_2, \ldots, a_j \in \{1, 2, \ldots i\} \text{ and } a_1 < a_2 < \ldots < a_j\} & \text{if } j \leq i \\ \phi & \text{if } j > i \end{cases}$

 - The composition of two morphisms $(a_1, a_2, \ldots, a_j)^i : i \to j$ and $(b_1, b_2, \ldots, b_k)^j : j \to k$ :

$$(b_1, b_2, \ldots, b_k)^j \circ (a_1, a_2, \ldots, a_j)^i = (a_{b_1}, a_{b_2}, \ldots, a_{b_k})^i.$$

Now, the Mobius function μ of $C_n$ is given by

$$\mu((a_1, a_2, \ldots, a_j)^i) = (-1)^{i-j}.$$

225